\documentclass[12pt]{amsart} 

\newtheorem{theorem}{Theorem}[section]
\newtheorem{lemma}[theorem]{Lemma}
\newtheorem{proposition}[theorem]{Proposition}
\newtheorem{corollary}[theorem]{Corollary}
\newtheorem{conjecture}[theorem]{Conjecture}

\newtheorem{example}[theorem]{Example}
\newtheorem{remark}[theorem]{Remark}
\newtheorem{question}[theorem]{Question}

\def\II{{\mathcal I}}
\def\OO{{\mathcal O}}

\def\mm{{\mathbb M}}
\def\CC{{\mathbb C}}
\def\PP{{\mathbb P}}
\def\NN{{\mathbb N}}

\newcommand{\Syz}{\mathrm{Syz}}

\begin{document}

\title[Implicitization in the presence of base points]{Implicitization
of surfaces in $\PP^3$ in the presence of base points}


\author{Laurent Bus\'e}
\address{UNSA, UMR 6621, Parc Valrose, BP 71, 06108 Nice Cedex 02, France}
\email{lbuse@unice.fr}

\author{David Cox}
\address{Department of Mathematics and Computer Science, Amherst
College, Amherst, MA 01002-5000, USA}
\email{dac@cs.amherst.edu}

\author{Carlos D'Andrea} 
\address{Miller Institute for Basic Research in Science and Department
of Mathematics, University of California at Berkeley, Berkeley, CA
94720-3840} 
\email{cdandrea@math.berkeley.edu}

\begin{abstract}
We show that the method of moving quadrics for implicitizing surfaces
in $\PP^3$ applies in certain cases where base points are present.
However, if the ideal defined by the parametrization is saturated,
then this method rarely applies.  Instead, we show that when the base
points are a local complete intersection, the implicit equation can be
computed as the resultant of the first syzygies.
\end{abstract}

\keywords{implicitization, syzygies, base points}

\subjclass{Primary 14Q10; Secondary 13D02, 68U07}

\maketitle

\section{Introduction}

Let $x(s,t,u),\, y(s,t,u),\,z(s,t,u)$ and $w(s,t,u)$ be homogeneous
polynomials of degree $n$ such that the parametrization
\begin{equation}
\label{uno}
X=\frac{x(s,t,u)}{w(s,t,u)},\
Y=\frac{y(s,t,u)}{w(s,t,u)},\ Z=\frac{z(s,t,u)}{w(s,t,u)}
\end{equation}
defines a surface in $\PP^3$.  The implicitization problem consists
in the computation of a homogeneous polynomial $P(X,Y,Z,W)$ whose
vanishing defines the projective closure of this surface.

The implicit equation can always be found using Gr\"obner bases. 
However, complexity issues mean that in practice, this method is
rarely used in geometric modeling, especially in situations where
real-time modeling is involved.  A more common method for finding the
implicit equation is to eliminate $s,t,u$ by computing the resultant
of the three polynomials
\[
x(s,t,u)-Xw(s,t,u),\ y(s,t,u)-Yw(s,t,u),\ z(s,t,u)-Zw(s,t,u).
\]
But in many applications, the resultant vanishes identically due to
the presence of base points, which are points $(s_0:t_0:u_0)\in\PP^2$
such that
\[
x(s_0,t_0,u_0)= y(s_0,t_0,u_0)= z(s_0,t_0,u_0)=w(s_0,t_0,u_0)=0
\]
(see \cite{CGZ} and the references therein).  

In \cite{SC}, Sederberg and Chen introduced a new technique for
finding the implicit equation \eqref{uno} called the \emph{method of
moving quadrics}.  This method is based in the construction of a
matrix $\mm$ whose entries are the coefficients in the monomial basis
of certain syzygies of the ideal $I = \langle x,y,z,w\rangle\subset
\CC[s,t,u]$, and syzygies of $I^2$.  The determinant of this matrix
is---under suitable assumptions---the implicit equation.  Having a
determinantal representation of the implicit equation is useful for
geometric modeling because of well-known algorithms for computing
symbolic determinants.  There is also considerable theoretical
interest in knowing when a resultant \cite{SZ, WZ, DD, DE} or an implicit
equation \cite{Be} can be represented as a single determinant.

Until now, the method of moving quadrics has been proved valid only in
the case where there are no base points (see \cite{CGZ,D,C}).  The
motivation for this paper is twofold:\ first, we wanted to prove the
validity of this method in the presence of base points under suitable
algebraic conditions on $I$, and second, we were curious what tools
and concepts from commutative algebra would be required.

Our results provide a positive answer to the first open question given
in the last section of \cite{CGZ}, in the sense that as the number of
base points of the parametrization increases, so does the number of
moving planes which occur in the matrix $\mm$.  Moreover, we show that
if the number of base points is greater than or equal to the degree of
the parametrization, then the implicit equation may be computed as the
determinant of a smaller matrix than the one proposed in \cite{CGZ}.

One can check also that our method, when applied to the case of no
base points, recovers the results of \cite{CGZ}.  Hence our method may
be regarded as a generalization of \cite{CGZ}.  When base points are 
present, our methods require that they be a local complete 
intersection.  The main theoretical tool used in the proof is the 
regularity of a homogeneous ideal.

In the second part of the paper, we turn to the case where the ideal
$I =\langle x,y,z,w\rangle$ is saturated.  We show in
Proposition~\ref{ch} that the method of moving quadrics works only if
the degree of the parametrization is $3$.  So other methods will be
needed.  Here, the key observation \cite{C} is that in the saturated
case, the syzygy module of $x,y,z,w$ is a free $\CC[s,t,u]$-module
with $3$ generators.  If we regard a syzygy $(A,B,C,D)$ as a
polynomial $AX+BY+CZ+DW$, then we show in Theorem~\ref{mt} that when
$I$ is a local complete intersection, we can recover the implicit
equation by taking the resultant to these syzygies.  This allows us to
regard a basis of the syzygy module as a generalization of the
$\mu$-basis for curves given in \cite{CSC} (see also \cite{CZS}).

In general, the search of formulas for implicitization rational
surfaces with base points is a very active area of research due to the
fact that, in practical industrial design, base points show up quite
frequently.  In \cite{MC}, a perturbation is applied to resultants in
order to obtain a nonzero multiple of the implicit equation.  In the
recent paper \cite{B}, a new projection operator called the
\emph{residual resultant} (introduced in \cite{BEM}) is developed for
computing the implicit equation when the base points locus is a local
complete intersection.  Recently, Abdallah Al-Amrani informed us that
the notes \cite{J} show how to compute the implicit equation as the
determinant of the approximation complexes discussed by Vasconcelos in
\cite{V}.

The paper is organized as follows.  In Section~\ref{dos}, we present
and discuss the formal structure of the method of moving quadrics in
terms of syzygies.  We show that there are only two possible sizes for
the matrix of moving planes and moving quadrics.  Then, in
Section~\ref{extension} we prove that under suitable assumptions on
$I$, the method actually computes the implicit equation.  We
illustrate our results with some examples.

In Section~\ref{saturated}, we discuss the case where the ideal is a
saturated local complete intersection.  We prove Theorem~\ref{mt},
which asserts that the resultant of a basis of the syzygy module gives
the implicit equation raised to a power equal to the degree of the
parametrization.

The paper concludes with some open questions in Section~\ref{fc}
suggested by the results of Sections~\ref{extension}
and~\ref{saturated}.  Appendices A and B give technical results about
basepoints and regularity used in Section~\ref{extension}.

\section{The formal structure of the method}
\label{dos}
Let $R:=\CC[s,t,u]$.  A syzygy on $I = \langle x,y,z,w\rangle \subset
R$ is a linear form
\[
aX+bY+cZ+dW\in R[X,Y,Z,W]
\]
such that $ax+by+cz+dw=0$.  This is a \emph{moving plane}.  In the
same way, we will define a syzygy on $I^2$ as a quadratic form in
$R[X,Y,Z,W]$ such that it vanishes when the variables are substituted
by the polynomials $x,y,z,w$.  These syzygies are called \emph{moving
quadrics}.

The method, as described in \cite{CGZ} for the case where the projective 
variety 
$V(I)$ is empty, consists in fixing a degree $d$ (in that case
$d=n-1$) and constructing a matrix $\mm$ of size
$\big(\genfrac{}{}{0pt}{}{d+2}{2}\big)$---the number of monomials in
three variables of degree $d$---having in their rows the coefficients
in the monomial basis of a basis of the syzygies of degree $d$ in the
variables $s,t,u$ on $I$ and---if there remains space---some linearly
independent syzygies on $I^2$, also having degree $d$ in $s,t,u$.

It is straightforward to verify that $\det(\mm)$ is a homogeneous
polynomial in $\CC[X,Y,Z,W]$ such that it vanishes on the surface
\eqref{uno}.  If this polynomial is not identically zero, then it must
be a multiple of $P(X,Y,Z,W)$.  Under appropriate assumptions on the
parametrization, one can show that this determinant gives a non-zero
constant times $P(X,Y,Z,W)$ (see \cite{CGZ,C}).

Suppose that $V(I)$ consists of only finitely many points (possibly the empty
set).  Assume also that the parametrization \eqref{uno} is proper and
that $V(I)$ is a local complete intersection (possibly empty).  Then,
it is well-known that the degree of $P(X,Y,Z,W)$ is equal to
$n^2-\deg(V(I))$, where $\deg(V(I)) = \dim_\CC(R/I)_k$ for $k\gg0$
(see \cite{C}).

We want to find moving planes and moving quadrics of degree $d$ such
that the determinant of the above matrix $\mm$ equals the implicit
equation of the surface.  To see what conditions the degree $d$ must
satisfy, consider the exact sequence of $\CC$-vector spaces:
\begin{equation}
\label{grad}
0\rightarrow \Syz(I)_{d}\rightarrow R_{d}^4\stackrel{B}{\rightarrow}
R_{d+n} \rightarrow \left(R/I\right)_{d+n}\rightarrow 0.
\end{equation}
Here, $\Syz(I)_{d}$ is the $\CC$-vector space of all syzygies of
degree $d$ on $I$ in the variables $s,t,u,$ and the map $B$ is
given by $(x,y,z,w).$

Let $m:=\dim \Syz(I)_{d}$ and $i:=\dim (R/I)_{d+n}$.  Since
\eqref{grad} is exact, we obtain
\begin{equation}
\label{eq1}
i + 4\,\begin{pmatrix} d+2\\ 2\end{pmatrix} = m +
\begin{pmatrix} d+n+2\\2\end{pmatrix}.
\end{equation}
Denote by $\mm$ the matrix of moving planes and moving quadrics 
as explained at the
beginning of this section. It is of size
$\big(\genfrac{}{}{0pt}{}{d+2}{2}\big)$, where $m$ of the rows are 
homogeneous of degree one in $X,Y,Z,W$ and the remaining are of 
degree $2$.  We want the determinant of $\mm$ to equal $P(X,Y,Z,W)$
(up to a nonzero constant).  Comparing degrees, we get the equation
\begin{equation}
\label{eq2}
m + 2\left(\begin{pmatrix} d+2\\2\end{pmatrix}-m\right)=
n^2-\deg(V(I)),
\end{equation}
with the additional condition
\begin{equation}
\label{ineq}
\begin{pmatrix} d+2\\2\end{pmatrix}-m\geq0,
\end{equation}
which says that the number of syzygies of degree $d$ on $I$ is
less than or equal to the size of $\mm$.

However, if we compare $\deg(V(I)) = \dim_\CC (R/I)_k$ for $k\gg0$
with $i = \dim_\CC (R/I)_{d+n}$, it makes sense to also assume that
\begin{equation}
\label{degdn}
i = \dim_\CC(R/I)_{d+n} = \deg(V(I)).
\end{equation}
Combining this with \eqref{eq2} gives
\begin{equation}
\label{eq3}
m + 2\left(\begin{pmatrix} d+2\\2\end{pmatrix}-m\right)=
n^2-i,
\end{equation}
and solving equations \eqref{eq1} and \eqref{eq3} in $d$ and $m$ leads
to the following solutions:
\begin{equation}
\label{dsols}
\begin{aligned}
d=n-1,\ &\text{in which case $m=n+i$,}\\
d=n-2,\ &\text{in which case $m=i-n$.}
\end{aligned}
\end{equation}
{}From this we see that, in the case where there are few base points,
the only possibility is $d=n-1$.

\begin{remark}
{\rm \cite{CGZ} treats the case when the parametrization has no base
points.  As just noted, this implies $d = n-1$.  Furthermore,
$\deg(V(I)) = 0$ in this case, so that \eqref{degdn} is equivalent to
the surjectivity of the map $B$ in \eqref{grad}.  In \cite{CGZ}, the
matrix of $B$ is denoted $MP$, so that the surjectivity of $B$ means
that $MP$ has maximal rank.  This is a part of the hypothesis of
Theorem~5.2 of \cite{CGZ}.}
\end{remark}

\section{Extension of the Method}
\label{extension}
In this section, we will extend the method of moving quadrics to the
case where base points are present.  In order to do this, we impose
the following \emph{base point conditions} on the input polynomials:
\begin{enumerate}
\item[BP1:] $x(s,t,u),\,y(s,t,u),z(s,t,u)$ and $w(s,t,u)$ are
homogeneous of degree $n$ and linearly independent over $\CC$.
\item[BP2:] $V(I)$ consists of a finite number of points and 
$\deg(V(I))$ equals the sum of the multiplicities of the distinct 
points in the locus $V(I)$.
\item[BP3:] There is $d\in\{n-2,n-1\}$ such that $\dim_\CC(R/I)_{d+n} =
\deg(V(I))$.
\item[BP4:] $w\in \mbox{sat}(x,y,z)$ (where ``$\mbox{sat}$'' denotes
saturation).
\item[BP5:] $\dim {\Syz(x,y,z)}_d=0$, where $d$ is as in BP3.
\end{enumerate}
We can explain these conditions as follows:
\smallskip

1.\ Condition BP1 is obvious, except possibly for the linear
independence.  For this, observe that a linear relation among
$x,y,z,w$ implies that the image of the parametrization is a plane.
This case is trivial.
\smallskip

2.\ The finiteness of $V(I)$ in Condition BP2 is equivalent to
assuming that $x,y,z,w$ have no common factors.  Also, the degree
formula for the image of the parametrization given in \cite{C}
involves the sum of the multiplicities of the base points.  
Finally note that $\deg(V(I))$ equals the sum of the multiplicities 
of the distinct points in the locus $V(I)$ if and only if $V(I)$ is a local complete intersection.
\smallskip

3.\ Condition BP3 is explained by \eqref{degdn} and \eqref{dsols} from
Section~\ref{dos}.  The surprise is that BP3 is equivalent to the
following regularity condition:

\smallskip

\begin{enumerate}
\item[BP3${}'$:] There is $d\in\{n-2,n-1\}$ such that $I$ is
$(d+n)$-regular.
\end{enumerate}
\smallskip

\noindent This follows from Theorem~\ref{apBthm2} since
$d\in\{n-2,n-1\}$ implies that $d+n \ge 2n-2$.  Regularity will play
an important role in the proof of Theorem~\ref{mpmq}.  See \cite{BS}
for a discussion of regularity.
\smallskip

4.\ Since $V(I)$ is a local complete intersection,
Corollary~\ref{apAcor} implies that we can obtain condition BP4 by
replacing the input polynomials with generic linear combinations of
them.
\smallskip

5.\ Consider the exact sequence
\[
0\rightarrow \Syz(x,y,z)_d \rightarrow \Syz(x,y,z,w)_d \rightarrow R_d,
\]
where the first map sends $(A,B,C)$ to $(A,B,C,0)$ and the second
sends $(A,B,C,D)$ to $D$.  This shows that BP5 implies the inequality
\eqref{ineq} for the given value of $d$.  Also note that in the case
where there are no base points, BP5 is satisfied provided that the
homogeneous resultant of $x,y,z$ is different than zero.  This is
shown in \cite[Lemma~5.1]{CGZ}.
\smallskip

In order to see the independence between the conditions, consider the
following examples.

\begin{example}{\rm
Take $x=s^5$, $y=t^5$, $z=su^4$, and $w=st^2u^2$.  Here, we have $n=5$
and $V(I)$ is the local complete intersection consisting of the point
$(0:0:1)$ of multiplicity $5$.  Thus Conditions BP1 and BP2 are
satisfied.  However, $d = 3$ or $4$ implies $d+n=8$ or $9$, yet one
can check easily with {\tt Macaulay 2} that regularity of $I$ is $10$.
Thus BP3${}'$ and its equivalent BP3 fail in this case.}
\end{example}

\begin{example}
\label{f5}{\rm
Set $x=su^2$, $y=t^2(s+u)$, $z=st(s+u)$, and $w=tu(s+u)$.  For this
parametrization, we have $n=3$ and $V(I)$ is the local complete
intersection consisting of the three points $(1:0:0)$, $(0:1:0)$ and
$(0:0:1)$ of respective multiplicity $2,3$ and $1$.  The implicit
equation is hence a cubic surface.  One can compute with {\tt Macaulay
2} that $I$ is $2n-2=4$-regular.  It follows easily that Conditions
BP1--BP4 are satisfied by taking $d=1$.  However, one can prove that
$\dim_{\CC}(\Syz(x',y',z')_1)=1$, where $x',y',z'$ are generic linear
combinations of $x,y,z,w$.  Thus BP5 does not hold in this example,
even if we use generic linear combinations of $x,y,z,w$.  We will see
later that the method of moving quadrics fails in this case.}
\end{example}

\subsection{Construction of the moving plane coefficient matrix}
\label{mpcm}
Let $x,y,z,w$ satisfy BP1--BP5 and consider the following algorithm:
\[
\begin{array}{l}
I_w:=\emptyset; \\
\Omega:= \{ s^it^ju^{d-i-j}(x,y,z),\ 0\leq i+j\leq d\};\\
\Gamma_w:=\{s^it^ju^{d-i-j}w,\ 0\leq i+j\leq d\}.
\end{array}
\]
While $\Gamma_w\neq\emptyset$
\begin{itemize}
\item Select a column $s^it^ju^{d-i-j}w$ from $\Gamma_w,$ and remove
it from $\Gamma_w;$
\item If $s^it^ju^{d-i-j}w$ is linearly independent from the columns
in $\Omega,$ then add it to $\Omega;$ otherwise, add $(i,j)$ to $I_w.$
\end{itemize}

Observe that at the beginning of the algorithm, the set $\Omega$ is
linearly independent by Condition BP5.

\begin{remark}
\label{cuenta}{\rm
It is straightforward to check that at the end of the algorithm,
$|I_w|=m =$ the dimension of the syzygies of degree $d$ on
$x,y,z,w$.}
\end{remark}

Also note that in the case were there are no base points, this
algorithm constructs the matrix denoted by $MP_I$ in \cite{CGZ}.

Now define $MP_{Iw}$ to be the coefficient matrix of the polynomials
\[
\begin{array}{l}
 s^it^ju^{d-i-j}(x,y,z),\ 0\leq i+j\leq d,\\
s^it^ju^{d-i-j}w, \ (i,j)\notin I_w.
\end{array}
\]
Observe that $MP_{Iw}$ is a submatrix of $MP$ having the same rank as
$MP$ and is maximal with this property.  Thus $MP_{Iw}$ has maximal
rank.

\subsection{The moving quadrics coefficient matrix}
Let $MQ$ be the coefficient matrix of the polynomials
\[
s^it^ju^{d-i-j}(x^2,y^2,z^2,xy,xz,yz,xw,yw,zw,w^2),
\]
and let $MQ_{Iw}$ be the submatrix determined by the polynomials
\begin{equation}
\label{following}
\begin{array}{l}
 s^it^ju^{d-i-j}(x^2,y^2,z^2,xy,xz,yz),\ 0\leq i+j\leq d,\\
s^it^ju^{d-i-j}(xw, yw,zw)\ (i,j)\notin I_w.
\end{array}
\end{equation}
The Theorem 5.1 in \cite{CGZ} may be extended as follows.

\begin{theorem}
\label{mpmq}
Let $x,y,z,w$ satisfy BP1--BP5 and construct $MQ_{Iw}$ as above.  Then
$MQ_{Iw}$ has maximal rank.  Furthermore, the columns of $MQ_{Iw}$ are
a basis of the $\CC$-vector space $I^2_{d+2n}$.
\end{theorem}

\begin{proof}
Suppose that there exist homogeneous polynomials of degree $d$, say
$p_1(s,t,u),\dots,p_9(s,t,u)$, such that
\begin{equation}
\label{cuatro}
p_1x^2+p_2y^2+p_3z^2+p_4xy+p_5xz+p_6yz+p_7xw+p_8yw+p_9zw=0,
\end{equation}
where the exponents of the monomials which appear in $p_7,p_8,p_9$ are
of the form $(i,j,d-i-j), (i,j)\notin I_w$.

Rewrite equation \eqref{cuatro} as
\begin{equation}
\label{cinco}
(p_1x+p_4y+p_5z+p_7w)x+(p_2y+p_6z+p_8w)y+(p_3z+p_9w)z=0.
\end{equation}
Equation \eqref{cinco} is a syzygy on $x,y,z$.

Condition BP4 implies that $V(x,y,z) = V(I)$ as subschemes of $\PP^2$.
It follows that \eqref{cinco} is a syzygy which vanishes at the base
point locus $Z=V(x ,y,z)$ in the sense of \cite{CS} (i.e., a syzygy
$(a_1,a_2,a_3)$ on $x,y,z$ vanishes on $V(x,y,z)$ if $a_i\in
\mathrm{sat}(x,y,z)$ for all $i$).

As the ideal generated by $x,y,z$ is a local complete intersection by
Condition BP2, Theorem 1.7 of \cite{CS} implies that these syzygies
are ``Koszul syzygies'' in the sense of \cite{C,CS}. Thus, there are
polynomials $a,b,c$ of degree $d$ such that
\[
\begin{aligned}
p_1x+p_4y+p_5z+p_7w &= ay+bz,\\
p_2y+p_6z+p_8w&= -ax+cz,\\
p_3z+p_9w&= -bx-cy.
\end{aligned}
\]
Since the exponents of the monomials in $p_9$ are not in $I_w,$ the
third equality tells us that the columns of $MP_{Iw}$ are linearly
dependent unless $p_3=p_9=b=c=0,$ but if this happens, then the second
equality will be $p_2y+p_6z+p_8w=-ax$.  Since $MP_{Iw}$ has maximal
rank, it follows that $a=p_2=p_6=p_8=0$.  But then the first equation
implies that $MP_{Iw}$ doesn't have maximal rank.  This contradiction
proves that $MQ_{Iw}$ has maximal rank.

It follows the columns of $MQ_{Iw}$ are $\CC$-linearly independent,
and they lie in the space of polynomials of degree $d+2n$ which belong
to $I^2$.  In order to see that they generate all of $I^2_{d+2n}$, we
argue as follows.

As $I^2 \subset \mathrm{sat}(I^2)$, we clearly have $I^2_{d+2n}\subset
\mathrm{sat}({I^2})_{d+2n}$.  We will show that the columns of
$MQ_{Iw}$ are actually a basis of $\mathrm{sat}(I^2)_{d+2n}$.

Since $V(I)$ consists of a finite number of points, and $I$ is
$(d+n)$-regular by BP3, Corollary 5 of \cite{Chl} implies that
$\mathrm{sat}(I^2)$ is $(d+2n)$-regular.  As in Section~\ref{dos}, we
let $i$ denote the dimension of $(R/I)_{d+n}$.  Using the exact
sequence $0\to I/I^2\to R/I^2\to R/I\to0$ and the fact that $I$ is a
local complete intersection of codimension two (see the proof of
Theorem 2.4 in \cite{CS}), one can show that the Hilbert polynomial of
$R/I^2$ is equal to $3i$.

As $I^2$ and $\mathrm{sat}(I^2)$ have the same Hilbert polynomial, and
as $\mathrm{sat}(I^2)$ is $(d+2n)$-regular, we see that
$\dim_{\CC}\big(R/\mathrm{sat}(I^2)\big)_{d+2n}=3i$, so the dimension
of $\mathrm{sat}(I^2)_{d+2n}$ is equal to
$\big(\genfrac{}{}{0pt}{}{d+2n+2}{2}\big) - 3i$.  Using
\eqref{following} and Remark~\ref{cuenta}, we see that the rank
of $MQ_{Iw}$ is equal to
\[
3\begin{pmatrix} d+n+2\\ 2\end{pmatrix}-3\begin{pmatrix} d+2\\
2\end{pmatrix}- 3i.
\]
This number equals $\dim_{\CC}(\mathrm{sat}(I^2)_{d+2n})$ when $d=n-2$
or $n-1$.  So the polynomials in \eqref{following} are actually a
basis of $\mathrm{sat}(I^2)_{d+2n}$.  From this the last part of the
theorem follows straightforwardly.
\end{proof}

\begin{remark}{\rm
The argument of Theorem~\ref{mpmq} shows that
\[
I^2_{d+2n} = \mathrm{sat}(I^2)_{d+2n}.
\]
By Theorem~\ref{apBthm2}, it follows that $I^2$ is $(d+2n)$-regular.}
\end{remark}

\subsection {The matrix of moving planes and quadrics}
\label{matriz}
As in \cite{CGZ}, we can obtain a square matrix $\mm$ of size
$\big(\genfrac{}{}{0pt}{}{d+2}{2}\big)$ whose rows  contain the
coefficients of $m$ linearly independent moving planes of degree $d$,
indexed by $(a,b)\in I_w$, and linearly independent moving quadrics
indexed by $(a,b)\notin I_w,\, 0\leq a,b,\,0\leq a+b\leq d$.

More precisely, the rows of $\mm$ are indexed by monomials of
degree $d$ in $3$ variables.  These rows are described as follows.  We
first construct $\big(\genfrac{}{}{0pt}{}{d+2}{2}\big)-m$ linearly
independent moving quadrics of the form
\[
Q_{i,j} := W^2 s^it^ju^{d-i-j} + \,\mbox{terms without}\, W^2,\
(i,j)\notin I_w.
\]
This is done by writing $s^it^ju^{d-i-j}w^2$ as a linear combination
of the polynomials in \eqref{following}.  We can do this since the
columns of $MQ_{Iw}$ generate $I^2_{d+2n}$ by Theorem~\ref{mpmq}.

To complete the matrix, we then find $m$ linearly independent moving
planes, which we write in the form
\[
P_{i,j}:= W s^it^ju^{d-i-j} + \,\mbox{terms not involving}\,
s^at^bu^{d-a-b}W,\,(a,b)\in I_w,
\]
for every $(i,j)\in I_w$.

The entries of the matrix $\mm$ are the coefficients of these moving
quadrics and moving planes.  By ordering its rows and columns
appropriately, we may assume that $\mm$ has the following form:
\[
\mm=\begin{pmatrix}
W^2+\cdots&&&&&\\
&\ddots&&&&\\
&&W^2+\cdots&&&\\
&&&W+\cdots&&\\
&&&&\ddots&\\
&&&&&W+\cdots\\
\end{pmatrix}.
\]
The first rows of $\mm$ consist of the coefficients of the moving
quadrics $Q_{i,j},\,(i,j)\notin I,$ and the last rows are the
coefficients of the moving planes $P_{i,j},\,(i,j)\in I$.

The following theorem is an extension of Theorem 5.2 in \cite{CGZ}.

\begin{theorem}
Suppose that $x,y,z,w$ satisfy BP1--BP5 and that the surface is
properly parametrized.  Then $\det(\mm)$ gives the implicit equation of
the parametric surface up to a nonzero constant.
\end{theorem}

\begin{proof}
The determinant has total degree
\[
2\left(\begin{pmatrix}d+2\\2\end{pmatrix}-m\right) + m,
\]
which by \eqref{eq2} equals
\[
n^2-i.
\]
This is the degree of the implicit equation since the parametrization
is generically one-to-one.

Checking the diagonal of $\mm$, we can see that the determinant of
$\mm$ has the term $W^{n^2-i}$, provided that it does not cancel with
other term of the same form.  But it is straightforward to see that
this is the highest power of $W$ which appears in the expansion of the
determinant.  So, $\det(\mm)\neq0$, and it is easy to see that
it vanishes whenever the point $(X,Y,Z,W)$ lies on the parametric surface
because each row represents a moving plane or quadric that follows the
surface.

It now follows easily that $\det(\mm)$ is the implicit equation of
the surface.
\end{proof}

\begin{example}{\rm
Take $x=st$, $y=u^2$, $z=s^2+tu$, and $w=tu$.  Here, the implicit
equation is $W^2(Z-W)-X^2Y,$ and the zero locus of $x,y,z,w$ in $\PP^2$ is
$\{(0:1:0)\}$, which is a base point of multiplicity $1$.

All the base point conditions are satisfied here.  The degree of the
parametrization is $n=2$ and the degree of the implicit equation is
$3=2^2-1$. The unique value of $d$ possible here is $d=n-1=1$.  The
matrix $MP$ in the lexicographic order $s>t>u$ is
\[
MP=\begin{pmatrix}
0&1&0&0&0&0&0&0&0&0\\
0&0&0&1&0&0&0&0&0&0\\
0&0&0&0&1&0&0&0&0&0\\
0&0&0&0&0&1&0&0&0&0\\
0&0&0&0&0&0&0&0&1&0\\
0&0&0&0&0&0&0&0&0&1\\
1&0&0&0&1&0&0&0&0&0\\
0&1&0&0&0&0&0&1&0&0\\
0&0&1&0&1&0&0&0&0&0\\
0&0&0&0&1&0&0&0&0&0\\
0&0&0&0&0&0&0&1&0&0\\
0&0&0&0&0&0&0&0&1&0
\end{pmatrix}.
\]
The rows of this matrix correspond to the coefficients in the monomial
basis of the polynomials $(s,t,u)(x,y,z,w)$.  It is straightforward to
verify that the last three rows (corresponding to $(s,t,u)w$) are
linear combinations of the previous rows, so we may choose
$I_w=\{(1,0),(0,1),(0,0)\}$.  This gives
\[
MP_{Iw}=\begin{pmatrix}
0&1&0&0&0&0&0&0&0&0\\
0&0&0&1&0&0&0&0&0&0\\
0&0&0&0&1&0&0&0&0&0\\
0&0&0&0&0&1&0&0&0&0\\
0&0&0&0&0&0&0&0&1&0\\
0&0&0&0&0&0&0&0&0&1\\
1&0&0&0&1&0&0&0&0&0\\
0&1&0&0&0&0&0&1&0&0\\
0&0&1&0&1&0&0&0&0&0
\end{pmatrix},
\]
which has maximal rank. This matrix gives three linearly
independent moving planes of degree $1$:
\[
\begin{aligned}
P_{0,0}&=uW-tY,\\
P_{1,0}&=sW-uX,\\
P_{0,1}&=tW-tZ+sX.
\end{aligned}
\]
In this case, as $i=1$ and $(n^2-n)/2-i = 0$, so that there are no
moving quadrics to consider.  Then the matrix of moving planes and
moving quadrics is
\[
\mm = \begin{pmatrix}
W&0&-X\\
X&W-Z&0\\
0&Y&-W
\end{pmatrix},
\]
and the determinant of this matrix gives the implicit equation.}
\end{example}

\begin{example}{\rm
Let $x=s^3$, $y=t^2u$, $z=s^2t+u^3$, and $w=stu$.  One can check that
all conditions are satisfied and $V(I)=\{(0:1:0)\}$ with multiplicity
$2$.  Again, the only possibility is $d=n-1$.  With the aid of {\tt
Maple}, we found the following five moving planes of degree $2$:
\[
\begin{aligned}
P_{0,0}&=u^2W+t^2X-stZ,\\
P_{0,1}&=tuW-suY,\\
P_{1,1}&=stW-s^2Y,\\
P_{2,0}&=s^2W-tuX,\\
P_{0,2}&=t^2W-stY,
\end{aligned}
\]
and the moving quadric $Q_{1,0}=suW^2-u^2XY$.  This gives
\[
\mm=\begin{pmatrix}
W&0&0&0&-X&0\\
-Y&W&0&0&0&0\\
0&0&W^2&0&0&-XY\\
0&-Y&0&W&0&0\\
0&0&-Y&0&W&0\\
0&-Z&0&X&0&W
\end{pmatrix}.
\]
One computes that $\det(\mm)=W^7-X^2Y^3ZW+X^3Y^4$, which is the
implicit equation of the parametric surface.}
\end{example}

\begin{example}{\rm
We present here a case where $d=n-2$.  This example is taken from
\cite{SC}.  Consider the following parametrization of a cubic surface
with $6$ base points:
\[
\begin{aligned}
  x &= s^{2} {t}+2 t^{3}+s^{2} {u}+4 {s} {t} {u}+4 t^{2} {u}+3 {s}
  u^{2}+2 {t} u^{2}+2 u^{3}, \\
  y &= -s^{3}-2 {s} t^{2}-2 s^{2} {u}-{s} {t} {u}+{s} u^{2}-2 {t}
  u^{2}+2 u^{3},\\
  z &= -s^{3}-2 s^{2} {t}-3 {s} t^{2}-3 s^{2} {u}-3 {s} {t} {u}+2 t^{2}
  {u}-2 {s} u^{2}-2 {t} u^{2},\\
  w &= s^{3}+s^{2} {t}+t^{3}+s^{2} {u}+t^{2} {u}-{s} u^{2}-{t} u^{2}-u^{3}.
\end{aligned}
\]
One can check with {\tt Macaulay 2} that $I$ is saturated, local
complete intersection and its regularity is $3$.  As shown in
\cite{SC}, we have the following basis of syzygies of degree
$d=n-2=1:$
\begin{equation}
\label{ant}
\begin{array}{l}
sX+tY+uZ,\\
s(Y+W)+t(2Y-Z)+u(Y+2W),\\
s(Z-Y)+t(-X+2W)+u(X-Y).
\end{array}
\end{equation}
The first syzygy shows that $\Syz(x,y,z)_d\ne0$, so that Condition
BP5 is not verified.  But if we consider $x,y,w$ instead, then it is
straightforward to check that all conditions are satisfied, and the
method produces the following matrix of moving planes (again there are
no moving quadrics to consider here):
\[
\mm = \begin{pmatrix}
Z-Y&-X+2W&X-Y\\
-Y-W&Z-2Y&-Y-2W\\
X&Y&Z
\end{pmatrix}.
\]
The determinant of this matrix is the determinant of the
matrix of syzygies in \eqref{ant}, which has been shown in \cite{SC}
to be the implicit equation of the surface.}
\end{example}

\begin{example}{\rm
In this example we focus on Condition BP5 and show that the method of
moving surfaces may fail if $\Syz(x',y',z')_d$ is nonzero when $x',
y', z'$ are generic linear combinations of $x,y,z,w$.  We take the
parametrization given in Example~\ref{f5}, and $d=1$.  The following
is a basis of syzygies of degree $1$:
\begin{equation*}
\begin{array}{l}
sW-uZ,\\
tW-uY,\\
tZ-sY.
\end{array}
\end{equation*}
Following the method of moving quadrics, this gives the $3\times 3$
matrix $\mm$ whose rows are given by these syzygies.  However, one
easily sees that $\det(\mm)$ is identically $0$ in this case and hence
it is not an implicit equation.}
\end{example}

\begin{remark}{\rm More generally, suppose that $x,y,z,w$ is a
parametrization which satisfies BP1--BP4 but fails BP5 even after a
generic coordinate change, as in the previous example.  Then it is
easy to see that the method of moving quadrics must fail in one of two
ways.  To see this, recall that BP5 implies the equality
\[
\begin{pmatrix}d+2\\2\end{pmatrix} - m \ge 0.
\]
So when BP5 fails, this inequality may fail, which means that the
number of linearly independent moving planes is greater than the
number of rows of $\mm$.  But even when the above inequality
holds, there are still problems, which we explain as follows. When
we replace $x,y,z,w$ with generic linear combinations, the
implicit equation of the surface must contain $W^{n^2-i}$.
However, since $\Syz(x',y',z')_d$ is nonzero, it follows that at
least one row of $\mm$ will not contain $W$, so that $W$ appears
to the power at most $n^2-i-1$ in $\det(\mm)$.  This contradiction
shows that the method of moving quadrics fails.}
\end{remark}

\section{The saturated case}
\label{saturated}

Now we will concentrate on the case where $I$ is a saturated local
complete intersection, with $V(I)$ consisting in a finite number of
points.  We will show that the method of moving quadrics rarely
applies and that when it does fail, it can often be replaced with a
nice resultant.

In this situation, it is well-known (see \cite[Prop.~5.2]{C}) that
$R/I$ is Cohen-Macaulay and the syzygy module $\Syz(I)$ is a free
graded $\CC[s,t,u]$-module.  We also have the following resolution of
$I$ (see \cite{C}):
\begin{equation}
\label{complex}
0\!\rightarrow\! R(-n-\mu_1)\oplus R(-n-\mu_2)\oplus R(-n-\mu_3)
\!\stackrel{A}{\rightarrow}\! R(-n)^4\!\stackrel{B}{\rightarrow}\!
I\!\rightarrow\! 0,
\end{equation}
where $\mu_1+\mu_2+\mu_3=n,$ the map $B$ is given by $(x,y,z,w),$ and
the columns of $A$ give three syzygies of degrees $\mu_1,\mu_2,\mu_3$
respectively which are free generators of $\Syz(I).$

\subsection{Limitations on the method of moving quadrics}
The following proposition shows that in the saturated case, the method
described in the previous section can be used only for very low
degrees.

\begin{proposition}
\label{ch}
If $I$ is saturated and satisfies Conditions BP1--BP5, then the method
of moving quadrics works only for $d=n-2$ and $n\leq 3$.
\end{proposition}

\begin{proof}
Let us first prove that the method does not work for $d =n-1$.  By
\eqref{dsols}, this implies $m=n+i$, and we also have $i =
\frac12(n^2+\mu_1^2+\mu_2^2+\mu_3^2)$ by \cite[Proposition 5.3]{C}.
Then the inequality \eqref{ineq} becomes
\[
\frac{(n+1)n}{2}\geq n+\frac12(n^2+\mu_1^2+\mu_2^2+\mu_3^2),
\]
which is impossible for positive values of $n$.

Now consider the case $d=n-2$. Here, $m = i-n$, and we have the same
formula for $i$.  Thus the inequality \eqref{ineq} becomes
\[
\frac{n(n-1)}{2}\geq \frac12(n^2+\mu_1^2+\mu_2^2+\mu_3^2)-n,
\]
which is equivalent to $n\geq \mu_1^2+\mu_2^2+\mu_3^2$.  As
$\mu_1+\mu_2+\mu_3=n$, we obtain
\[
n\geq \mu_1^2+\mu_2^2+\mu_3^2\geq \frac{n^2}{3}.
\]
This shows that $n$ must be at most $3$.

{}From here, it is now easy to see that the only nontrivial case is
$n=3$ and $\mu_1=\mu_2=\mu_3=1$ (otherwise, the surface will be a
plane).  In this case, we have $m=3$, which is the number of monomials
of degree $1$ in $3$ variables, and the matrix $\mm$ is the matrix of
the basis of syzygies on $(x,y,z,w)$ of degree one.  The determinant
of this matrix gives the implicit equation. This can be proved by
hand, or regarded as a special case of Theorem~\ref{mt} (see
Corollary~\ref{cor}).
\end{proof}

\subsection{The implicit equation as a resultant} In \cite{CSC}, an
exact sequence similar to \eqref{complex} was used to represent the
implicit equation of a parametric curve as the resultant of the
homogeneous polynomials which were free generators of the syzygy
module (see also \cite{C} for an exposition of this).  We will discuss
whether these results extend to surfaces in the saturated case.

Consider again the exact sequence (\ref{complex}).
Write
\[
A=\begin{pmatrix}
p_1&q_1&r_1\\
p_2&q_2&r_2\\
p_3&q_3&r_3\\
p_4&q_4&r_4\\
\end{pmatrix}.
\]
This means that the polynomials
\[
\begin{aligned}
{\bf p}&= p_1 X + p_2 Y + p_3 Z + p_4 W,\\
{\bf q}&= q_1 X + q_2 Y + q_3 Z + q_4 W,\\
{\bf r}&= r_1 X + r_2 Y + r_3 Z + r_4 W
\end{aligned}
\]
are syzygies of degrees $\mu_1,\mu_2,\mu_3$ in the variables $s,t,u$,
which generate the syzygy module of $x,y,z,w$.  Let ${\rm
Res}_{\mu_1,\mu_2,\mu_3}(\cdot,\cdot,\cdot)$ be the homogeneous
resultant of three homogeneous polynomials of degrees $\mu_1, \mu_2,
\mu_3$ as defined in \cite{CLO}.  One may ask whether ${\rm
Res}_{\mu_1,\mu_2,\mu_3}({\bf p},{\bf q},{\bf r})$ computes a power of
the implicit equation, as in the case of curves. Unfortunately, the
following example shows that this is not always the case.

\begin{example}
\label{vanonempty}{\rm
Let
\[
\begin{aligned}
x&=st^3  - s^4  - 2s^2t^2  + s^2tu+4s^3 t - 2 t^3u,\\
y&= s^2tu- s^3t- 2s^3u+3st^2u-t^3u, \\
z&=s^3u-st^3- 4s^2tu+6t^3u-st^2u, \\
w&=s^3u-3st^3-2st^2u+6s^2t^2+ t^4 -ts^3.
\end{aligned}
\]
The ideal generated by these polynomials is saturated, and $A$
is the following matrix:
\[
\begin{pmatrix} s&2s&tu\\
s&t&s^2\\
2t&s&t^2\\
t&3t&su
\end{pmatrix}.
\]
All entries in this matrix vanish under the substitution
$s\mapsto0,\,t\mapsto0$, so the homogeneous resultant of the first
syzygies will be identically zero due to the fact that the polynomials
${\bf p},{\bf q},{\bf r}$ have the common root $(0:0:1)$ in projective
space.}
\end{example}

However, the ideal of Example \ref{vanonempty} is not a local complete
intersection.  If we add this hypothesis (which is part of Condition
BP2 from Section~\ref{extension}), then we get the following nice
result.

\begin{theorem}
\label{mt}
Assume that $I$ is saturated and satisfies Conditions BP1 and BP2.  Then
\begin{equation}
\label{res}
{\rm Res}_{\mu_1,\mu_2,\mu_3}({\bf p},{\bf q},{\bf r}) =
{P(X,Y,Z,W)}^{h},
\end{equation}
where $h$ is the degree of the parametrization and $P(X,Y,Z,W) = 0$ is
an implicit equation of the surface.
\end{theorem}

\begin{proof}
Let $S = V(P) \subset \PP^{3}$ be the Zariski closure of the image of
the parametrization \eqref{uno}.  To prove the theorem, first suppose
that the resultant vanishes at a point $(X_0:Y_0:Z_0:W_0)\in\PP^3$.
This means that the system of equations in variables $s,t,u$ given by
\begin{equation}
\label{system}
\begin{aligned}
p_1X_0+p_2Y_0+p_3Z_0+p_4W_0&=0,\\
q_1X_0+q_2Y_0+q_3Z_0+q_4W_0&=0,\\
r_1X_0+r_2Y_0+r_3Z_0+r_4W_0&=0
\end{aligned}
\end{equation}
has a non-trivial solution $(s_0,t_0,u_0)$.  We will show that
\begin{equation}
\label{toprove}
(X_0:Y_0:Z_0:W_0)\in S \cup {\textstyle\bigcup_{p \in V(I)}} L_p,
\end{equation}
where $L_p$ is a line.  Since the right-hand side is a proper
subvariety of $\PP^3$, this will prove that ${\rm
Res}_{\mu_1,\mu_2,\mu_3}({\bf p},{\bf q},{\bf r})$ is a nonzero
polynomial and hence has zero locus of pure codimension $1$.  Since
$L_p$ has codimension $2$ and $V(I)$ is finite, this will prove that
the zero locus lies in $S$.

Given the solution $(s_0,t_0,u_0)$ of \eqref{system}, we can
specialize the variables $s,t,u$ to $s_0,t_0,u_0$ in the exact
sequence \eqref{complex}.  This transforms \eqref{complex} into a
complex of vector spaces
\[
0\rightarrow\CC^3\stackrel{A_0}{\rightarrow}
\CC^4\stackrel{B_0}{\rightarrow} \CC\rightarrow0,
\]
where $A_0$ is the matrix $A$ specialized, and $B_0 =
(X_0,Y_0,Z_0,W_0)$.  As $B_0$ is surjective and---because of
(\ref{complex})---$B_0A_0=0$, we see that the complex is exact if and
only if $A_0$ is injective.  By the Hilbert-Burch Theorem, the maximal
minors of $A_0$ are $x(s_0,t_0,u_0)$, $y(s_0,t_0,u_0)$,
$z(s_0,t_0,u_0)$, $w(s_0,t_0,u_0)$.  So, if we are outside of the zero
locus of $I$, the complex is exact, and the determinant of the complex
is non-zero (see \cite[Appendix A]{GKZ} for a definition of the
determinant of a complex).  Moreover, applying the Cayley formula for
computing this determinant with respect to the monomial bases, we get
the following:
\[
\begin{aligned}
x(s_0,t_0,u_0)&=X_0 D,\\
y(s_0,t_0,u_0)&=Y_0 D,\\
z(s_0,t_0,u_0)&=Z_0 D,\\
w(s_0,t_0,u_0)&=W_0 D,
\end{aligned}
\]
where $D$ is the determinant of the complex.  From here, it is easy to
see that the point $(X_0:Y_0:Z_0:W_0)$ belongs to the surface $S$.

However, if $p = (s_{0}:t_{0}:u_{0}) \in V(I)$, then the above
argument fails.  To see what happens in this case, we first study the
rank of the specialized matrix $A_{0}$.  Localizing \eqref{complex} at
$p$ gives
\begin{equation}
\label{localcx}
0 \rightarrow \mathcal{O}_{p}^3 \stackrel{A}{\rightarrow}
\mathcal{O}_{p}^{4}
 \stackrel{B}{\rightarrow} \mathcal{I}_{p} \rightarrow 0,
\end{equation}
where $\mathcal{O} = \mathcal{O}_{\PP^2}$ and $\mathcal{I}$ is the
ideal sheaf associated to $I$.  Since $\mathcal{I}_{p} \subset
\mathcal{O}_{p}$ is a complete intersection, the minimal resolution of
$\mathcal{I}_{p}$ is of the form:
\[
0 \rightarrow \mathcal{O}_{p} \rightarrow \mathcal{O}_{p}^{2}
\rightarrow \mathcal{I}_{p} \rightarrow 0.
\]
This means that \eqref{localcx} is isomorphic to the exact sequence
obtained from the minimal resolution by adding the trivial complex
\[
0 \rightarrow \mathcal{O}_{p}^2 =
\mathcal{O}_{p}^{2}
 \rightarrow 0 \rightarrow 0.
\]
Hence $A$ has a $2\times2$ minor which doesn't vanish at $p$.  In
other words, the matrix $A_{0}$ has rank $\ge 2$.  It follows that
substituting $p$ into \eqref{system} gives a system of linear equations
of rank $\ge 2$ when regarded as equations in $X_{0},Y_0, Z_0, W_0$.
However, we also know that the rank is $< 3$ since the $3\times 3$
minor of $A$ are $x,y,z,w$, which vanish at $p \in V(I)$.
Projectively, this means that $(X_0:Y_0:Z_0:W_0)$ belongs to the line
$L_p \subset \PP^3$ defined by substituting $p$ into \eqref{system}.
This completes the proof of \eqref{toprove}.

The next step is to show that the resultant vanishes on $V(P)$, and
for this, it is enough to show that it vanishes on a Zariski dense
subset.  For instance, we can take the image of the parametrization
\[
\begin{aligned}
\PP^2\setminus V(I) &\to \PP^3,\\
(s:t:u) &\mapsto (x(s,t,u):y(s,t,u):z(s,t,u):w(s,t,u)).
\end{aligned}
\]
For $(X_0:Y_0:Z_0:W_0)$ in the image, we can find
$(s_0:t_0:u_0)\in\PP^2$ in the preimage.  It is straightforward to
check that the syzygies ${\bf p},\,{\bf q}, \,{\bf r}$ vanish after
the specialization of all the variables.  Thus ${\rm
Res}_{\mu_1,\mu_2,\mu_3}({\bf p},{\bf q},{\bf r})$ vanishes at
$(X_0:Y_0:Z_0:W_0)$.

Since $P$ is irreducible, it follows that
\[
{\rm Res}_{\mu_1,\mu_2,\mu_3}({\bf p},{\bf q},{\bf
r})=c{P(X,Y,Z,W)}^{\delta}
\]
for some $\delta\in\NN$ and a non-zero constant multiplier $c$.  
To see that $\delta$ is the degree of the
parametrization, note that by \cite{C}, the degree of the surface is
equal to $\mu_1\mu_2 + \mu_1\mu_3 + \mu_2\mu_3$, which is equal to the
degree of ${\rm Res}_{\mu_1,\mu_2,\mu_3}({\bf p},{\bf q},{\bf r})$ in
the variables $X,Y,Z,W$.  By the degree formula (see
\cite[Appendix]{C}), this number must be $h$ times the degree of
$P(X,Y,Z,W)$.
\end{proof}

\begin{corollary}
\label{cor}
If $n=3$ and $\mu_1=\mu_2=\mu_3=1$, then the implicit equation is the
determinant of the first syzygy module.
\end{corollary}

\section{Open Questions}
\label{fc}

\begin{question}{\rm Most of the base point conditions imposed on the
ideal $I$ in Section~\ref{extension} were needed in order to prove
that matrix $\mm$ has nonzero determinant.  A straightforward
computation shows that---for the degrees $d$ of 
Section~\ref{dos}---there is a natural map
\begin{equation}
\label{mapi}
\begin{aligned}
{\Syz(I)_{d}}^4 &\to \Syz(I^2)_{d},\\
(S_1,S_2,S_3,S_4)&\mapsto S_1X+S_2Y+S_3Z+S_4W
\end{aligned}
\end{equation}
whose cokernel has dimension greater than or equal to
$\big(\genfrac{}{}{0pt}{}{d+1}{2}\big)-m$.  Thus, if inequality
\eqref{ineq} holds, then we can fill $\mm$ with moving quadrics which
do not come from the previous map.  It is easy to see that, in order
to have $\det(\mm)\neq0$, the moving quadrics of $\mm$ must not belong 
to the image of \eqref{mapi}.  Is this a sufficient
condition?  This would make it easier to compute the implicit equation,
and we would have a general result with fewer conditions on the base
points.}
\end{question}

\begin{question}{\rm In order to construct matrix $\mm$ we used
\emph{all} moving planes of a given degree.  Can we make this
condition weaker, i.e., can we use matrices which use some but not all
moving planes of degree $d$?}
\end{question}

\begin{question}{\rm In the situation of Theorem~\ref{mt}, one can ask
how the resultant ${\rm Res}_{\mu_1,\mu_2,\mu_3}({\bf p},{\bf q},{\bf
r})$ relates to the implicit equation $P = 0$ when $I$ is not
necessarily a local complete intersection.  In general, one can show
that if $h$ is the degree of the parametrization, then
\begin{equation}
\label{genhpres}
\begin{aligned}
h \deg(P) &= \mu_1\mu_2 + \mu_1\mu_3 + \mu_2\mu_3 - {\textstyle
\sum_{p \in V(I)}} (e_p-d_p),\\
        &=\deg\big({\rm Res}_{\mu_1,\mu_2,\mu_3}({\bf p},{\bf q},{\bf
r})\big)- {\textstyle
\sum_{p \in V(I)}} (e_p-d_p),
\end{aligned}
\end{equation}
where
\begin{align*}
e_p &= \text{\rm multiplicity of}\ \II_p \subset \OO_p,\\
d_p &= \text{\rm degree of}\ \II_p \subset \OO_p.
\end{align*}
Here, $\II_p$ is the ideal of the local ring $\OO_p$ induced by $I$.

To analyze this, let $A$ be as in \eqref{complex} and let $V_i(A)
\subset \PP^2$ be the subscheme defined by the vanishing of the
$i\times i$ minors of $A$.  Then
\[
V_1(A) \subset V_2(A) \subset V_3(A) = V(I),
\]
where the last equality holds by the Hilbert-Burch Theorem.  Hence
there are three cases to consider:

\smallskip

\emph{Case 1:} $V_1(A) \ne \emptyset$.  In this situation, it is easy
to see that the resultant vanishes identically.  This is what happened
in Example~\ref{vanonempty}.

\smallskip

\emph{Case 2:} $V_2(A) = \emptyset$.  When $I$ has a resolution of the form
\eqref{complex}, it is easy to show that
\[
V_2(A) = \emptyset \Leftrightarrow I\ \text{is a local complete
intersection}.
\]
Hence this case is covered by Theorem~\ref{mt}.

\smallskip

\emph{Case 3:} $V_1(A) = \emptyset$ and $V_2(A) \ne \emptyset$.  When
$p \in V_2(A)$ is substituted into \eqref{system}, the resulting
system of linear equations has rank $1$ and hence defines a plane $H_p
\subset \PP^3$.  Then the argument used to prove \eqref{toprove} can
be modified so show that
\begin{equation}
\label{v2a}
V\big({\rm Res}_{\mu_1,\mu_2,\mu_3}({\bf p},{\bf q},{\bf r})\big) = S
\cup {\textstyle \bigcup_{p \in V_2(A)}} H_p.
\end{equation}
It follows that the resultant has extraneous factors in this case.  If
$\ell_p = 0$ is the equation of the plane $H_p$, then we have the
following conjectural formula for the resultant.

\begin{conjecture}
\label{conj}
Let $I \subset \CC[s,t,u]$ be generated by $x,y,z,w$ of degree $n$
such that $V(I)$ is finite.  Also assume that:
\begin{enumerate}
\item $I$ is saturated with free resolution given by \eqref{complex}.
\item $V_1(A) = \emptyset$ and $V_2(A) \ne \emptyset$.
\end{enumerate}
Then, up to a nonzero constant, we have
\begin{equation}
\label{conjres}
{\rm Res}_{\mu_1,\mu_2,\mu_3}({\bf p},{\bf q},{\bf r}) = P(X,Y,Z,W)^h\,
{\textstyle \prod_{p \in V_2(A)}}\, \ell_p^{e_p-d_p}.
\end{equation}
\end{conjecture}

This conjecture is compatible with \eqref{v2a} since $p \in V_2(A)
\Rightarrow \II_p$ is not a complete intersection $\Rightarrow e_p >
d_p$.  Furthermore, \eqref{genhpres} shows that each side of
\eqref{conjres} has the same degree since $e_p = d_p$ for $p \in V(I)
\setminus V_2(A)$.

For $p \in V_2(A)$, notice that $e_p - d_p$ measures how far $V(I)$ is
{}from being a complete intersection at $p$.  Hence,
Conjecture~\ref{conj}, if true, would show that resultants are
sensitive to subtle features of the base point locus.  Also, how does
$e_p - d_p$ relate to the subscheme structure of $V_2(A) \subset
\PP^2$ at $p$?

Finally, suppose that $V_1(A)$ is nonempty (as in Case 1) and is a
local complete intersection.  Is there a version of
Conjecture~\ref{conj} which uses the residual resultant (see
\cite{B})?

The following example illustrates how extraneous components can arise
as predicted by Conjecture \ref{conj}.}

\begin{example} {\rm
Consider the parametrization given by
\begin{align*}
  x &= {s} t^{2}-t^{3}-{t} u^{2}, \\
  y &= t^{3}-{s} {t} {u}-t^{2} {u}+{t} u^{2}+u^{3}, \\
  z &= {s} {t} {u}-2 {t} u^{2},\\
  w &= t^{2} {u}-2 {t} u^{2}+u^{3}.
\end{align*}
Here, we have $n=3$.  Using {\tt Macaulay2}, it is easy to
compute that $I=\langle x,y,z,w\rangle$ is saturated with free
resolution \eqref{complex} where the matrix $A$ is given by:
\[
A=\begin{pmatrix}
  {-u}&{t}& 0\\
  {-u}& -s+{t}+{u}&{-u}\\
  {t}-u&-s-t+2 u&{-u}\\
  {u}&{s}-u&{t}+{u}\\
\end{pmatrix}.
\]
One can also show that $V(I)$ consists of points $p=(1:0:0) $ and
$q=(2:1:1)$.  The point $p$ is in $V_2(A)$, the point $q$ is in $V(I)
\setminus V_2(A)$, and one easily checks that $V_1(A)$ is empty.
Hence the hypothesis of Conjecture \ref{conj} is satisfied.

Since $\mu_1 = \mu_2 = \mu_3 = 1$, the degree of $I$ is
$\frac12(n^2+\mu_1^2+\mu_2^2+\mu_1^3) = 6$ by \cite{C}.  Thus $d_p =
d_q = 1$.  The resultant of ${\bf p},{\bf q},{\bf r}$ is
the determinant
\[
\det\begin{pmatrix}
  0&-Y-Z+{W}& 0\\
  {Z}&{X}+{Y}-Z&{W}\\
  -X-Y-Z+{W}&{Y}+2 Z-W&-Y-Z+{W}\\
\end{pmatrix},
\]
which can be factored into
\[
(-{Y} {Z}-Z^{2}+{X} {W}+{Y} {W}+2 {Z} {W}-W^{2}) ({Y}+{Z}-W).
\]
One can check that the first factor gives the implicit equation, which
hence has degree $2$, and that the parametrization has degree $1$. 
{}From this, we know that $3^2-e_p-e_q=2$.  Moreover, since $q \in V(I)
\setminus V_2(A)$, $\II_q$ is a local complete intersection and hence
$d_q=e_q$.  Furthermore, since $p \in V_2(A)$, $\II_q$ is not a
local complete intersection and hence $e_p>d_p$.  This implies that
$e_p-d_p=1$, so that Conjecture~\ref{conj} predicts that the resultant
has a single extraneous component of multiplicity $1$.  This is
confirmed by the above factorization.}
\end{example}
\end{question}

\begin{question}{\rm The proof of Theorem~\ref{mt} also shows that
each base point $p \in V(I)$ blows up to a line $L_p \subset S$.
What happens if we drop the hypothesis that $V(I)$ is a local complete
intersection from the theorem?  For example, suppose that $V_1(A) =
\emptyset$ and $p \in V_2(A)$.  Then the resultant has an extraneous
factor $\ell_p$ as in Conjecture~\ref{conj}.  If $H_p \subset \PP^3$
is the plane defined by $\ell_p = 0$, then presumably $p \in V_2(A)$
blows up to the plane curve $S \cap H_p \subset H_p$.  But then what
happens if $p \in V_1(A)$?  Does $p$ blow up to a space curve which
doesn't lie any plane?  All of this indicates an interesting relation
between the geometry of a parametrization and the structure of
various subschemes of its base point locus.}
\end{question}

\section*{Acknowledgments}

The authors are grateful to Abdallah Al-Amrani and Marc Chardin for
giving us some intuition about regularity.  We also thank Haohao Wang
for pointing out a missing hypothesis in the original version of Lemma
B.2. 

Part of this work was done while the third author was a Postdoctoral
Fellow at the Institut National de Recherche en Informatique et en
Automatique (INRIA) in Sophia-Antipolis, France, partially supported
by Action A00E02 of the ECOS-SeTCIP French-Argentina bilateral
collaboration.

\appendix

\section{A theorem about basepoints\\ }\label{apA}

We begin with the following general result.

\begin{theorem}
\label{apAthm}
Let $X$ be a Cohen-Macaulay variety of dimension $d$ and let
$\mathcal{L}$ be a line bundle on $X$.  Also assume that $L \subset
H^{0}(X,\mathcal{L})$ is a subspace such that ${V}(L) \subset
X$ is a 0-dimensional subscheme.  If $s_{0},\dots,s_d \in L$ are
generic, then:
\begin{enumerate}
\item ${V}(s_{0},\dots,s_d) = {V}(L)$ as sets.
\item ${V}(s_{0},\dots,s_d)$ and ${V}(L)$ have the
same multiplicity at all points.
\end{enumerate}
Furthermore, if ${V}(L)$ is a local complete intersection, then we
have ${V}(s_{0},\dots,s_d) = {V}(L)$ as subschemes.
\end{theorem}

\begin{proof}
Let $m = \dim_{\CC}(L)$.  If $m \le d+1$, then $s_{0},\dots,s_d$ span
$L$ and it follows that ${V}(s_{0},\dots,s_d) = {V}(L)$ as subschemes
of $X$.  Hence we may assume that $m > d+1$.

{}From $L$ we get the morphism $\varphi: X \setminus {V}(L) \to
\PP(L^{*}) \simeq \PP^{m-1}$.  The image of $\varphi$ is a
constructible set of dimension at most $d$.  Hence we can find a
linear subvariety $P$ of codimension $d+1$ which is disjoint from the
image.  We can write $P$ as the intersection of $d+1$ generic
hyperplanes.  However, hyperplanes in $\PP(L^{*})$ are defined by
elements of $L$.  Thus $P$ is defined by $s_{0},\dots,s_d \in L$.
Furthermore, $P$ being disjoint from the image of $\varphi$ implies
that $s_{0},\dots,s_d$ don't vanish simultaneously on
$X\setminus{V}(L)$, i.e., ${V}(s_{0},\dots,s_d) \subset {V}(L)$.  The
other inclusion is obvious, which completes the proof of part (a) of
the theorem.

For part (b), fix $p \in {V}(L)$ and let $I_{p} \subset \OO_{p}$ be
the ideal generated by $L$ in the local ring $\OO_{p}$of $X$ at $p$.
Corollary 4.5.10 of Bruns and Herzog \cite{bh} implies that $\OO_{p}$
has a system of parameters which generates a reduction ideal $J_{p}$
for $I_{p}$.  Note that this system of parameters is a regular
sequence since $\OO_p$ is Cohen-Macaulay (Theorem 2.12 of \cite{bh}).
Furthermore, the proofs of Proposition 4.5.8, Theorem 1.5.17 and
Proposition 1.5.12 of \cite{bh} show that the system of parameters can
be chosen to be generic linear combinations of generators of $I_p$.
Since we can use a basis of $L$ as generators of $I_p$, it follows
that the system of parameters can be chosen to be generic elements of
$L$.  This system has $d$ elements since $\OO_p$ has dimension $d$.

It follows that $s_{0},\dots,s_{d-1}$ can be assumed to be a regular
sequence which generates a reduction ideal $J_{p}$ for $I_{p}$.
Furthermore, since this is true for generic elements of $L$ and
${V}(L)$ is finite, we can assume that this holds for all $p
\in {V}(L)$.

Now let $\tilde I_{p}$ be the ideal of $\OO_{p}$ generated by
$s_{0},\dots,s_d$.  Then we have the obvious inclusions
\[
J_{p} \subset \tilde I_{p} \subset I_{p},
\]
which gives the inequalities
\[
e(J_{p}) \ge e(\tilde I_{p}) \ge e(I_{p}).
\]
However, the first and third terms are equal since $J_{p}$ is a
reduction ideal for $I_{p}$.  This proves the desired equality of
multiplicities.

Finally, if $I_{p}$ is a complete intersection, then it
coincides with all of its reduction ideals (this is easy to prove).
Thus $J_{p} = I_{p}$, which by the above inclusions implies $\tilde
I_{p} = I_{p}$.   This shows that ${V}(s_{0},\dots,s_d)$
and ${V}(L)$ have the same scheme structure at $p$.  When
${V}(L)$ is a local complete intersection, this is true for
all of its points, and it follows that ${V}(s_{0},\dots,s_d)
= {V}(L)$ as schemes.
\end{proof}

\begin{corollary}
\label{apAcor}
Suppose $x,y,z,w \in \CC[s,t,u]$ are homogeneous of degree
$n$ with no common factor.  If we replace $x,y,z$ with generic
linear combinations of $x,y,z,w$, then
\begin{enumerate}
\item ${V}(x,y,z) = {V}(x,y,z,w)$ as sets.
\item ${V}(x,y,z)$ and ${V}(x,y,z,w)$ have the same
multiplicity at each point.
\end{enumerate}
Furthermore, if ${V}(x,y,z,w)$ is a local complete
intersection, then we have ${V}(x,y,z) = {V}(x,y,z,w)$ as
subschemes and $w \in \mathrm{sat}(x,y,z)$.
\end{corollary}

\begin{proof} Apply Theorem~\ref{apAthm} to $L =
\mathrm{Span}(x,y,z,w) \subset H^{0}(\PP^{2},
\mathcal{O}_{\PP^{2}}(n))$.  For the final assertion,
note that ${V}(x,y,z) = {V}(x,y,z,w)$ as subschemes $\Leftrightarrow
\mathrm{sat}(x,y,z) = \mathrm{sat}(x,y,z,w) \Leftrightarrow w \in
\mathrm{sat}(x,y,z)$.
\end{proof}

\section{A theorem about regularity \\}\label{apB}

We begin with a lemma in two variables.

\begin{lemma}
\label{apBlem1}
Let $\bar{I} \subset \CC[s,t]$ have $r$ minimal homogeneous 
generators of degree
$n$.  If $V(\bar{I}) = \emptyset$ in $\PP^1$, then $\bar{I}$ is
$m$-regular for all $m \ge 2n-r+1$.
\end{lemma}

\begin{proof}
Let $S = \CC[s,t]$.  The Hilbert syzygy theorem, together with a
Hilbert polynomial calculation and the fact that $\bar{I}_k = S_k$ for
$k \gg 0$, imply that $\bar{I}$ has a minimal graded free resolution
\[
0 \to \bigoplus_{i=1}^{r-1} S(-n-\mu_i) \to S(-n)^r \to \bar{I} \to 0,
\]
where $\mu_i \ge 1$ for all $i$ and $\sum_{i=1}^{r-1}
\mu_i = n$.

Since $S(-n)$ and $S(-n-\mu_i)$ have generators of degrees $n$ and
$n+\mu_i$ respectively, Definition~3.2(c) of \cite{bm} implies that
$\bar{I}$ is $m$-regular whenever $m \ge \max\{n,n+\mu_i-1\} = n +
\max\{\mu_i\} - 1$.  However, $\mu_i \ge 1$ and $\sum_{i=1}^{r-1}
\mu_i = n$ imply that for each $i$, we have $n \ge \mu_i + r - 2$.
This implies $2n - r + 1 \ge n + \max\{\mu_i\} - 1$, and the lemma
follows.
\end{proof}

For the rest of this appendix, we will study the regularity of certain
homogeneous ideals $I \subset R = \CC[s,t,u]$ using the inductive 
method found on
page 34 of \cite{bm}.  We begin with the following result.

\begin{lemma}
\label{apBlem2}
Let $I \subset R = \CC[s,t,u]$ have $r \ge 4$ minimal homogeneous
generators, all of degree $n$, and assume that $V(I) \subset \PP^2$ is
finite and the rational map from $\mathbb{P}^2$ to $\mathbb{P}^{r-1}$
given by the minimal generators is generically finite.  Given a
generic element of $\ell \in R_1$, let $I_\ell$ be the image of $I$ in
the quotient ring $R/\langle \ell\rangle$.  Then $I_\ell$ has at least
$3$ minimal generators.
\end{lemma}

\begin{proof}
Let $p_1,\dots,p_r$ be minimal homogeneous generators of $I$, 
where each $p_i$ has
degree $n$.  Then let $Z \subset \PP^{r-1} \times \PP(R_1) = \PP^{r-1}
\times \PP^2$ be defined by
\[
Z = \{([a_1,\dots,a_r],[\ell]) : \ell | a_1p_1 + \dots + a_rp_r\},
\]
and let $\pi_1 : Z \to \PP^{r-1}$ and $\pi_2 : Z \to \PP^2$ be the
natural projections.  Note that our hypothesis implies $n > 1$.

Since $V(I)\subset \PP^2$ is finite, we know that $p_1,\dots,p_r$ have
no common factors.  Thus the linear system of divisors given by
$a_1p_1 + \dots + a_rp_r = 0$ is reduced in the sense of \cite[p.\
130]{ii}, and the image of the rational map it determines has
dimension $2$ by hypothesis.  It follows by the Bertini theorem
\cite[Thm.\ 7.19]{ii} that the general member of the linear system is
irreducible.  Since $n > 1$, we conclude that $\pi_1^{-1}(p) =
\emptyset$ for a general point $\text{\bfseries\emph{a}} \in
\PP^{r-1}$.  Furthermore, if $\pi_1^{-1}(\text{\bfseries\emph{a}}) \ne
\emptyset$, then $\pi_1^{-1}(\text{\bfseries\emph{a}})$ is finite
since $a_1p_1 + \dots + a_rp_r$ is divisible by at most $n$ linear
forms.  Standard arguments then imply that $Z$ has dimension $\le
r-2$.

Now consider a generic $\ell \in \PP^2$ and let $(p_i)_\ell$ denote
the image of $p_i$ in $R/\langle \ell\rangle$.  If $\pi_2^{-1}(\ell) =
\emptyset$, then we are done since the $(p_i)_\ell$
are linearly independent in this case.  On the other hand, suppose
that $\pi_2^{-1}(\ell) \ne \emptyset$ when $\ell$ is generic.  Since
$Z$ has dimension $\le r-2$, it follows that
\[
\pi_2^{-1}(\ell)\ \text{has dimension at most}\  r-4
\]
for generic $\ell$.  However,
\[
\pi_2^{-1}(\ell) = \PP(\text{space of linear relations among the}\
(p_i)_\ell),
\]
so that for generic $\ell$, the space of linear relations among the
$(p_i)_\ell$ has dimension $\le r-3$.  This implies that at least $3$ of
$(p_i)_\ell$ are linearly independent for generic $\ell$.
\end{proof}

We now state our first main result.

\begin{theorem}
\label{apBthm1}
Let $I \subset R = \CC[s,t,u]$ have $r \ge 4$ minimal homogeneous
generators, all of degree $n$, and assume that $V(I) \subset \PP^2$ is
finite and the rational map from $\mathbb{P}^2$ to $\mathbb{P}^{r-1}$
given by the minimal generators is generically finite.  If
$\mathcal{I}$ is the associated sheaf on $\PP^2$, then:
\begin{enumerate}
\item $H^2(\mathcal{I}(k)) = 0$ for all $k \ge 0$.
\item $H^1(\mathcal{I}(k)) = 0$ for all $k \ge 2n-3$.
\end{enumerate}
\end{theorem}

\begin{proof}
Let $Z = V(I) \subset \PP^2$.  The statement for $H^2(\mathcal{I}(k))$
is then a trivial consequence of
\begin{equation}
\label{zex}
0 \to \II \to \OO_{\PP^2} \to \OO_Z \to 0
\end{equation}
and the vanishing of the higher cohomology of $\OO_Z$.

To prove the second statement, let $\ell$ be a generic element of
$R_1$.  Since $Z$ is finite, we may assume that $Z \cap V(\ell) =
\emptyset$.  By Lemma~\ref{apBlem2}, we may also assume that
\[
\bar{I} = I_\ell \subset S = R/\langle \ell\rangle
\]
has at least $3$ minimal homogeneous generators, all of degree $n$.  Then
Lemma~\ref{apBlem1} implies that $\bar{I}$ is $m$-regular for $m \ge
2n -3 + 1 = 2n-2$.  If $\bar{\II}$ is the sheaf associated to
$\bar{I}$, then by Definition 3.2(b) of \cite{bm}, we have
\begin{align}
\label{bm1}
{}&\text{$\bar{I}_k \to H^0(\bar{\II}(k))$ is an
isomorphism for $k \ge 2n-2$},\\ \label{bm2}
{}& H^1(\bar{\II}(k)) = 0\ \text{for}\ k \ge 2n-3.
\end{align}

We now use the argument of \cite[p.\ 34]{bm}.  Tensoring
\[
0 \to \OO_{\PP^2}(-1) \to \OO_{\PP^2} \to \OO_{\PP^1} \to 0
\]
with $\II(k)$ gives the exact sequence
\[
Tor_1^{\OO_{\PP^2}}(\II(k),\OO_{\PP^1}) \to \II(k-1) \to \II(k) \to
\bar{\II}(k) \to 0.
\]
However, $Tor_1^{\OO_{\PP^2}}(\II(k),\OO_{\PP^1})$ is supported on $V(\ell)
\simeq \PP^1$, and $\II(k)$ is locally free on $\PP^2 - V(I)$.  Then
$Tor_1^{\OO_{\PP^2}}(\II(k),\OO_{\PP^1}) = 0$ follows from $V(I) \cap
V(\ell) = \emptyset$.

\newcommand{\hz}{H^0(\II}
\newcommand{\da}{\downarrow}
\newcommand{\sa}{\!\to\!}

Thus we have an exact sequence
\[
0 \to \II(k-1) \to \II(k) \to
\bar{\II}(k) \to 0,
\]
whose long exact sequence in cohomology gives the commutative diagram
\[
\begin{array}{cccccccccccc}
         &\!\!&I_k    &\sa& \bar{I}_k       &\sa&0            &\\
         &\!\!&\da    &   &\downarrow       &   &\da          &\\
\hz(k-1))&\sa &\hz(k))&\sa&H^0(\bar{\II}(k))&\sa&H^1(\II(k-1))&\sa\\
 &&&&&&\\
&\to&H^1(\II(k))&\to&H^1(\bar{\II}(k))&\to&
\end{array}
\]
with exact rows.

Now suppose that $k \ge 2n-2$.  Then \eqref{bm2} implies that
$H^1(\bar{\II}(k)) = 0$.  Furthermore, $I_k \to \bar{I_k}$ is onto and
$\bar{I_k} \to H^0(\bar{\II}(k))$ is an isomorphism when $k \ge 2n-2$
by \eqref{bm1}.  Thus $H^0(\II(k)) \to H^0(\bar{\II}(k))$ is onto when
$k \ge 2n-2$.  Hence the above diagram gives an isomorphism
\[
H^1(\II(k-1)) \simeq H^1(\II(k)),\quad k \ge 2n-2.
\]
But we also know that $H^1(\II(k)) = 0$ for $k \gg 0$.  It
follows easily that
\[
H^1(\II(k-1)) = 0,\quad k \ge 2n-2.
\]
This implies the second statement of the theorem.
\end{proof}

Our second main result now follows easily.  Given $I$ as above, recall
that for $Z = V(I) \subset \PP^2$, we have
\[
\deg(Z) = \mathrm{dim}\, H^0(\OO_Z) = \mathrm{dim} (R/I)_k, \quad k \gg
0.
\]

\begin{theorem}
\label{apBthm2}
Let $I \subset R = \CC[s,t,u]$ have $r \ge 4$ minimal homogeneous
generators, all of degree $n$, and assume that $V(I) \subset \PP^2$ is
finite and the rational map from $\mathbb{P}^2$ to $\mathbb{P}^{r-1}$
given by the minimal generators is generically finite.  If $m \ge
2n-2$, then $I$ is $m$-regular if and only if $\mathrm{dim} (R/I)_m =
\deg(Z)$.
\end{theorem}

\begin{proof}
Observe that $m \ge 2n-2$ and Theorem~\ref{apBthm1} imply
$H^1(\II(m)) = 0$.  Hence \eqref{zex} gives the exact sequence
\begin{equation}
\label{mex}
0 \to H^0(\II(m)) \to R_m \to H^0(\OO_Z) \to 0, \quad m \ge 2n-2.
\end{equation}
Now suppose that $I$ is $m$-regular.  This implies $I_m =
H^0(\II(m))$, and then $\mathrm{dim}\, (R/I)_m = \deg(Z)$ follows
easily from \eqref{mex}.

Conversely, suppose that $\mathrm{dim} (R/I)_m = \deg(Z)$.  Then $m
\ge 2n-2$ and Theorem~\ref{apBthm1} imply that $H^1(\II(m-1)) =
H^2(\II(m-2)) = 0$ (note that $m-2 \ge 0$ since $n \ge 2$).  By
Definition 3.2(a) of \cite{bm}, $I$ will be $m$-regular once we prove
that $I_m \to H^0(\II(m))$ is an isomorphism.  Furthermore, since this
map is injective, it suffices to show $\mathrm{dim}\, I_m =
\mathrm{dim}\, H^0(\II(m))$.  However, \eqref{mex} implies that
\[
\mathrm{dim}\, H^0(\II(m)) = \mathrm{dim}\, R_m - \deg(Z).
\]
Combining this with $\mathrm{dim} (R/I)_m = \deg(Z)$ immediately
implies that $\mathrm{dim}\, H^0(\II(m)) = \mathrm{dim}\, I_m$, and
the theorem is proved.
\end{proof}


\begin{thebibliography}{XXX}

\bibitem[Be]{Be} A.~Beauville.
\newblock {\em Determinantal hypersurfaces.\/}
\newblock Michigan Math.\ J.\ {\bf 48} (2000), 39--64.

\bibitem[Bu]{B} L.~Bus\'e.
\newblock {\em Residual resultant over the projective plane and the
implicitization problem.\/}
\newblock {\sl Proceedings of ISSAC 2001,\/} 48--55.

\bibitem[BEM]{BEM}
L.~Bus\'e, M.~Elkadi and B.~Mourrain.
\newblock{\em Resultant over the residual of a complete intersection.\/}
\newblock {\sl Effective Methods in Algebraic Geometry (Bath,
2000).\/}
J.\ Pure Appl.\ Algebra {\bf 164} (2001), 35--57.

\bibitem[BH]{bh}  W.\ Bruns and J.\ Herzog.
\newblock {\sl Cohen-Macaulay Rings.\/}
\newblock Cambridge U.\ Press, Cambridge, 1993.

\bibitem[BM]{bm} D.\ Bayer and D.\ Mumford.
\newblock {\em What can be computed in algebraic geometry?}
\newblock {\sl Computational Algebraic Geometry and Commutative
Algebra.\/}
\newblock Cambridge Univ.\ Press, Cambridge, 1993, 1--48.

\bibitem[BS]{BS} D. Bayer and M. Stillman.
\newblock {\em A criterion for detecting $m$-regularity.\/}
\newblock Invent.\ Math.\ {\bf 87} (1987), 1--11.

\bibitem[Ca]{Can} J.~F.~Canny.
\newblock {\em Generalized characteristic polynomials.\/}
\newblock J.\ Symbolic Comput.\ {\bf 9} (1990), 241--250.

\bibitem[Ch]{Chl} K.~A.~Chandler.
\newblock {\em Regularity of the powers of an ideal.\/}
\newblock Comm.\ Algebra {\bf 25} (1997), 3773--3776.

\bibitem[Co]{C} D.~Cox.
\newblock {\em Equations of parametric curves and surfaces via
syzygies.\/}
\newblock {\sl Symbolic Computation:\ Solving Equations in Algebra,
Geometry and Engineering.\/}
\newblock Contemporary Mathematics, vol.\ 286, AMS, Providence, RI,
2001, 1--20.

\bibitem[CGZ]{CGZ} D.~Cox, R.~Goldman and M.~Zhang.
\newblock {\em On the validity of implicitization by moving quadrics
for rational surfaces with no base points.\/}
\newblock J. Symbolic Comput.\ {\bf 29} (2000), 419--440.

\bibitem[CLO]{CLO} D.~Cox, J.~Little and D.~O'Shea.
\newblock {\sl Using Algebraic Geometry}.
\newblock Springer-Verlag, New York Berlin Heidelberg, 1998.

\bibitem[CS]{CS} D.~Cox and H.~Schenck.
\newblock {\em Local complete intersections in $\PP^2$ and Koszul syzygies.\/}
\newblock Proc.\ Amer.\ Math.\ Soc.\ {\bf 131} (2003), 2007--2014.

\bibitem[CSC]{CSC} D.~Cox, T.~Sederberg and F.~Chen.
\newblock {\em The moving line ideal basis of planar rational curves.\/}
\newblock Comput.\ Aided Geom.\ Des.\ {\bf 15} (1998), 803--827.

\bibitem[CZS]{CZS} F.~Chen, J.~Zheng, and T.~W.~Sederberg.
\newblock {\em The $\mu$-basis of a rational ruled surface.}
\newblock  Comput.\ Aided Geom.\ Des.\ {\bf 18} (2001), 61--72.

\bibitem[D]{D} C.~D'Andrea.
\newblock {\em Resultants and moving surfaces.\/}
\newblock J. Symbolic Comput.\ {\bf 31} (2001), 585--602.

\bibitem[DD]{DD} C.~D'Andrea and A.~Dickenstein.
\newblock {\em Explicit formulas for the multivariate resultant.\/}
\newblock J.\ Pure Appl.\ Algebra {\bf 164} (2001), 59--86.

\bibitem[DE]{DE} A.~Dickenstein and I.~Emiris.
\newblock {\em Multihomogenous resultant matrices.\/}
\newblock Proceedings of the
2002 International Symposium on Symbolic and Algebraic Computation
(Lille),
\newblock ACM, New York, 2002. 

\bibitem[E]{E} D.~Eisenbud.
\newblock{\sl Commutative Algebra with a View Toward Algebraic Geometry.\/}
\newblock Springer-Verlag, New York Berlin Heidelberg, 1995.

\bibitem[GKZ]{GKZ} I.~Gelfand,  M.~Kapranov and A.~Zelevinsky.
\newblock {\sl Discriminants, Resultants and Multidimensional Determinants.\/}
\newblock Birkh\"{a}user, Boston, 1994.

\bibitem[I]{ii} S.\ Iitaka.
\newblock {\sl Algebraic Geometry\/}.
\newblock Springer-Verlag, New York Berlin Heidelberg, 1982.

\bibitem[J]{J}  J.~P.~Jouanolou.
\newblock  Course DEA, 2000--2001, University of Strasbourg.

\bibitem[MC]{MC}  D.~Manocha and J.~F.~Canny.
\newblock {\em Algorithms for implicitizing rational parametric surfaces.\/}
\newblock Comput.\ Aided Geom.\ Des.\ {\bf 9} (1992), 25--50.

\bibitem[SC]{SC}  T.~W..~Sederberg and F.~Chen.
\newblock{\em Implicitization using moving curves and surfaces.\/}
\newblock {\sl  Proceedings of SIGGRAPH 1995,\/} 301--308.

\bibitem[SZ]{SZ} B.~Sturmfels and A.~Zelevinsky.
\newblock {\em Multigraded resultants of Sylvester type.\/}
\newblock J.\ Algebra {\bf 163} (1994), 115--127.

\bibitem[V]{V}  W.~V.~Vasconcelos
\newblock {\sl Arithmetic of Blowup Algebras.\/}
\newblock London Mathematical Society Lecture Note Series, vol.\ 195.
\newblock Cambridge University Press, Cambridge, 1994.

\bibitem[WZ]{WZ} J.~Weyman and A.~Zelevinsky. 
\newblock {\em Determinantal formulas for multigraded resultants.\/} 
\newblock J.\ Algebraic Geom.\ {\bf 3} (1994), 569--597.


\end{thebibliography}
\end{document}